\documentclass{amsart}

\usepackage{amsmath}  
\usepackage{amsfonts}
\usepackage{tikz}  
\usepackage{enumitem}

\usepackage{etoolbox}
\makeatletter
\patchcmd{\@maketitle}
  {\ifx\@empty\@dedicatory}
  {\ifx\@empty\@date \else {\vskip3ex \centering\footnotesize\@date\par\vskip1ex}\fi
   \ifx\@empty\@dedicatory}
  {}{}
\patchcmd{\@adminfootnotes}
  {\ifx\@empty\@date\else \@footnotetext{\@setdate}\fi}
  {}{}{}
\makeatother

\def\vec#1{\ensuremath{\mathchoice
    {\mbox{\boldmath$\displaystyle\mathbf{#1}$}}
    {\mbox{\boldmath$\textstyle\mathbf{#1}$}}
    {\mbox{\boldmath$\scriptstyle\mathbf{#1}$}}
    {\mbox{\boldmath$\scriptscriptstyle\mathbf{#1}$}}}}%
\def\char{\mbox{1\hspace{-.25em}l}}

\newcommand{\vect}[1]{\vec{#1}}
\def\tens#1{\relax\ifmmode\mathsf{#1}\else\textsf{#1}\fi}
\newcommand{\matr}[1]{\tens{#1}}

\newcommand{\Lip}{\mathrm{Lip}}

\newtheorem{theorem}{Theorem}[section]
\newtheorem{lemma}[theorem]{Lemma}
\newtheorem{proposition}[theorem]{Proposition}

\newtheorem{definition}[theorem]{Definition}

\theoremstyle{remark}
\newtheorem{remark}[theorem]{Remark}

\def\R{\hbox{\bf R}}
\def\Z{\hbox{\bf Z}}

\def\N{\hbox{\bf N}}

\newcommand{\ba}{\begin{eqnarray}}
\newcommand{\ea}{\end{eqnarray}}

\renewcommand{\N}{{\mathbb N}}
\renewcommand{\R}{{\mathbb R}}

\renewcommand{\Z}{{\mathbb Z}}

\newenvironment{Proofc}[1]{\smallskip\par\noindent\textsc{#1}\quad}%
  {\hfill$\Box$\bigskip\par}
\newenvironment{Proof}{\begin{Proofc}{Proof}}{\end{Proofc}}
\title[Gradient entropy estimate for diagonal hyperbolic systems]{Gradient entropy estimate\\ and convergence of a semi-explicit scheme\\
  for diagonal hyperbolic systems}

\author{L.~Monasse, R.~Monneau}
\address{CERMICS, Ecole des Ponts ParisTech, Universit\'e Paris-Est, F-77455 Marne-la-Vall\'ee, France}
\email{\{monassel,monneau\}@cermics.enpc.fr}

\date{\today}

\begin{document}

\begin{abstract}
In this paper, we consider diagonal hyperbolic systems with monotone continuous initial data.
We propose a natural semi-explicit and upwind first order scheme.
Under a certain non-negativity condition on the Jacobian matrix of the velocities of the system, 
there is a gradient entropy estimate for the hyperbolic system.
We show that our scheme enjoys a similar gradient entropy estimate at the discrete level.
This property allows us to prove the convergence of the scheme.
\end{abstract}

\keywords{Semi-explicit upwind scheme, diagonal hyperbolic systems, gradient entropy estimate, monotone initial data}

\maketitle

\section{Introduction}

%\subsection{Discussion of the origin of the problem}
In this paper, we are interested in diagonal hyperbolic systems with monotone continuous initial data,
and in their discretization.
In a first subsection, we present our framework for such hyperbolic systems. In a second subsection, 
we propose a natural semi-explicit scheme.
In a third subsection we give our main results, including the convergence of the scheme. In a fourth subsection, we recall the related literature.
Finally, in a fifth subsection, we give the organization of the paper.

\subsection{The continuous problem}

Let us consider the following diagonal hyperbolic system (in nonconservative form).
Let $\vect{u}:\mathbb{R}\times[0,T]\rightarrow\mathbb{R}^d$ be a
solution of:
\begin{equation}
  \frac{\partial u^\alpha}{\partial t} + \lambda^\alpha(\vect{u})\frac{\partial u^\alpha}{\partial
    x} = 0 \textrm{ in }\mathcal{D}'((0,+\infty)\times \R) \label{eqn:system}
\end{equation} 
with initial data
\begin{equation}\label{eq::2}
u^\alpha(0, \cdot)=u_0^\alpha,\quad \mbox{for}\quad 
\alpha=1,...,d
\end{equation}

In order to specify our conditions on the initial data, it will be useful to recall the definition of the Zygmund space:
$$L\log L(\R)=\left\{w\in L^1(\R),\quad \int_{\R} |w| \ln (e + |w|) < +\infty \right\}$$
which is a Banach space with the norm
\begin{equation*}\label{eq::12}
|w|_{L\log L(\R)}=\inf \left\{\mu>0,\quad \int_{\R} \frac{|w|}{\mu}\ln \left( e + \frac{|w|}{\mu}\right) \le 1\right\}
\end{equation*}
Then we will assume that the initial data satisfies
\begin{equation}\label{eq::r2}
\left\{\begin{array}{l}
u^\alpha_0 \quad \mbox{is bounded and non-decreasing},\\
(u^\alpha_0)_x \in L\log L(\R)
\end{array}\right| \quad \mbox{for}\quad \alpha=1,...,d
\end{equation}
In particular such initial data is continuous.

We equip from now on the vector space $\mathbb{R}^d$ with the
1-norm $|\vect{u}| = \sum_{\alpha=1}^d{|u^{\alpha}|}$. We assume that
\begin{equation}\label{eq::r1}
\vect{\lambda} \in C^1(\R;\R^d) \quad \mbox{with $\vect{\lambda}$ globally Lipschitz continuous}
\end{equation}
with a Lipschitz constant $\Lip(\vect{\lambda})$. 
In addition,  the symmetric part of the  jacobian matrix of $\vect{\lambda}$ is supposed to be non negative in the following sense:
\begin{equation}\label{eq::7}
\sum_{\alpha,\beta = 1,...,d}\xi_\alpha \xi_\beta \frac{\partial \lambda^\alpha}{\partial u^\beta}(u)  \ge 0 \quad \mbox{for every}\quad \vect{\xi}=(\xi_1,...,\xi_d) \in [0,+\infty)^d,\quad u\in \R^d
\end{equation}
where we notice that this inequality is required only for a subset of vectors $\vect{\xi} \in \R^d$ with non negative coordinates.
When $d=1$, this condition for Burgers type equations ensures that solutions associated to non decreasing and continuous initial data
will stay continuous for all positive times. Under assumption (\ref{eq::7}) for $d\ge 1$, we can recover a similar property:
 it is possible to show that the solutions formally satisfy the following inequality 
\begin{equation}\label{eq::r10}
  \frac{d}{dt}\int_{\mathbb{R}}{\sum_{\alpha=1}^{d}{u^\alpha_x\ln(u^\alpha_x)}}\leq 0
\end{equation}
Indeed, we refer the reader to  Theorem 1.1 and Remark 1.4 in \cite{EM1}, for  a precise statement.
In some cases (in particular to ensure the uniqueness of the solution), 
we will also assume that the system is strictly hyperbolic, i.e. $\lambda$ satisfies:
\begin{equation}\label{eq::8}
\lambda^\alpha(\vect{u}) < \lambda^{\alpha+1}(\vect{u}) \quad \mbox{for}\quad \alpha=1,...,d-1
\end{equation}
We also define the total variation of $\vect{u}$ at time $\tau$ on the open
interval $(a,b)$: 
\begin{equation*}
TV[\vect{u}(\tau);(a,b)] = \sup\left\{\sum_{\alpha=1,...,d}\int_a^b{-{u}^\alpha(\tau,x) \varphi^\alpha_x(x) dx}\right\}
\end{equation*}
where the supremum is taken over the set of functions $\varphi^\alpha \in C^1_c(a,b)$ satisfying  $|\varphi^\alpha(x)|\le 1$ for $x\in (a,b)$ and $ \alpha=1,...,d$. In the particular case where for each $\alpha=1,...,d$,  the function $u^\alpha(\tau,\cdot)$ belongs to $W^{1,1}_{loc}(\R)$ and is non
decreasing in space, then we simply have
$$TV[\vect{u}(\tau);(a,b)]  = \int_a^b{|\vect{u}_x(\tau,x)|dx}=|\vect{u}(\tau,b)-\vect{u}(\tau,a)|$$
In the special case where $(a,b)=\mathbb{R}$, we will simply write
$TV[\vect{u}(\tau)]$. 

\begin{definition}[\bf Continuous vanishing viscosity solutions]\label{defi::1}
  A function $\vect{u}\in [C([0,+\infty)\times\mathbb{R})]^d$ is a continuous vanishing viscosity
  solution of system (\ref{eqn:system})-(\ref{eq::2}) if $\vect{u}$ solves (\ref{eqn:system})-(\ref{eq::2}) and if the following integral estimate
  holds.

  There exist constants $C$, $\gamma$, $\eta>0$ such that, for every
  $\tau\geq 0$ and $a<\xi<b$, with $b-a\leq\eta$, one has the tame estimate
  \begin{equation}\label{eq:tameest}
    \limsup_{h\rightarrow 0^+}{\frac{1}{h}\int_{a+\gamma h}^{b-\gamma
        h}{|\vect{u}(\tau+h,x)-\vect{U}_{(\vect{u}(\tau);\tau,\xi)}(h,x)|dx}}\leq C\left(TV[\vect{u}(\tau);(a,b)]\right)^2 
  \end{equation}
  where $\vect{U}_{\vect{u}(\tau);\tau,\xi}$ is the solution of the
  linear hyperbolic Cauchy problem with frozen constant coefficients:
  \begin{equation*}
    \frac{\partial w^\alpha}{\partial t}
    +\lambda^{\alpha}(\vect{u}(\tau,\xi))\frac{\partial w^\alpha}{\partial x} =
    0,\quad\text{with }w^\alpha(0,x)=u^\alpha(\tau,x).
  \end{equation*}
\end{definition}

This definition is in El Hajj, Monneau \cite{EM2} and is an adaptation
of the definition of Bianchini, Bressan \cite{BB} 
(see also in the book of Dafermos \cite{D} and the tame
oscillation estimate for solutions constructed with the Front tracking
method; this last estimate is related but less precise than the tame
estimate (\ref{eq:tameest})).

We then recall the following result (see Theorem 1.1 and Remark 1.4 in \cite{EM1},\cite{EM2}).

\begin{theorem}{\bf (Existence, uniqueness)}\label{th::1}\\
Assume that initial data satisfies (\ref{eq::r2}), and that $\vect{\lambda}$ satisfies (\ref{eq::r1}) and (\ref{eq::7}).\\
\noindent {\bf i) (Existence)} Then there exists a function $\vect{u}\in \left(C([0,+\infty)\times \R)\right)^d$ with $\vect{u}_x\in \left(L^\infty((0,+\infty); L\log L (\R)\right)^d$, which is a continuous vanishing viscosity solution of (\ref{eqn:system})-(\ref{eq::2}) in the sense of Definition \ref{defi::1}.\\
\noindent {\bf ii) (Uniqueness)} If moreover the system is strictly hyperbolic, i.e. $\vect{\lambda}$ safisties (\ref{eq::8}), then there is uniqueness of the continuous vanishing viscosity solution $\vect{u}$ of (\ref{eqn:system})-(\ref{eq::2}) in the sense of definition \ref{defi::1}.
\end{theorem}

\subsection{The semi-explicit discretization}

To recover these properties on the discrete level, 
we consider a time step $\Delta t>0$ and a space step $\Delta x>0$ and consider $u^{\alpha,n}_i$ as an approximation of $u^\alpha(n\Delta t, i\Delta x)$. 
We propose the following semi-explicit discretization of the system:
\begin{equation*}\label{eq::r5}
  \forall \alpha\in\{1,\dots,d\},
  \left\{
  \begin{array}{r}
    \displaystyle
    \frac{u_i^{\alpha,n+1}-u_i^{\alpha,n}}{\Delta
      t}+\lambda^\alpha(\vect{u}_i^{n+1})\left(\frac{u_{i+1}^{\alpha,n}-u_{i}^{\alpha,n}}{\Delta
        x}\right) =0 \quad \text{if } \lambda^{\alpha}(\vect{u}_i^{n+1})\leq 0\\
    \displaystyle
    \frac{u_i^{\alpha,n+1}-u_i^{\alpha,n}}{\Delta
      t}+\lambda^\alpha(\vect{u}_i^{n+1})\left(\frac{u_{i}^{\alpha,n}-u_{i-1}^{\alpha,n}}{\Delta x}\right) =
    0 \quad \text{if } \lambda^\alpha(\vect{u}_i^{n+1})\geq 0
  \end{array}
  \right. 
\end{equation*}
It is a first-order upwind formulation, with the velocity $\vect{\lambda}(\vect{u})$
being implicit in time. We denote
$$\lambda_i^{\alpha,n+1}=\lambda^\alpha(\vect{u}_i^{n+1})$$ and we define
its positive and negative parts $(\lambda_i^{\alpha,n+1})_+$ and
$(\lambda_i^{\alpha,n+1})_-$ as follows:
\begin{equation*}
  (\lambda_i^{\alpha,n+1})_+ =
  \frac{1}{2}(\lambda_i^{\alpha,n+1}+|\lambda_i^{\alpha,n+1}|), \quad (\lambda_i^{\alpha,n+1})_- =
  \frac{1}{2}(|\lambda_i^{\alpha,n+1}|-\lambda_i^{\alpha,n+1})
\end{equation*}
Both $(\lambda_i^{\alpha,n+1})_+$ and
$(\lambda_i^{\alpha,n+1})_-$ are positive real numbers.
We can write the scheme in a more compact form:
\begin{equation}\label{eqn:schema}
  \frac{u_i^{\alpha,n+1}-u_i^{\alpha,n}}{\Delta
      t}-(\lambda_i^{\alpha,n+1})_-\left(\frac{u_{i+1}^{\alpha,n}-u_{i}^{\alpha,n}}{\Delta
        x}\right) + (\lambda_i^{\alpha,n+1})_+\left(\frac{u_{i}^{\alpha,n}-u_{i-1}^{\alpha,n}}{\Delta
        x}\right) =0
\end{equation}

In the sequel, we set:
\begin{equation}\label{eq::r9}
  \theta_{i+\frac{1}{2}}^{\alpha,n} = \frac{u_{i+1}^{\alpha,n}-u_i^{\alpha,n}}{\Delta x}
\end{equation}
which  is a discrete equivalent of
$u_x^{\alpha}$. 

For a fixed index $i_0$ and $N\in\mathbb{N}$, we denote 
\begin{equation}\label{eq::r11}
I_N(i_0) = \{i_0-N,\dots,i_0+N\}, 
\end{equation}
and we define $TV[\vect{u}^n;I_N(i_0)]$ the
total variation of $\vect{u}^n$ on the set of indices $I_N(i_0)$:
\begin{equation*}
  TV[\vect{u}^{n};I_N(i_0)] = \sum_{\alpha=1}^{d}{\sum_{i\in I_N(i_0)}{|u_{i+1}^{\alpha,n}-u_{i}^{\alpha,n}|}}.
\end{equation*}
The total variation of $\vect{u}^n$ on $\mathbb{Z}$ is simply noted
$TV(\vect{u}^n)$.

\subsection{Main results}

%\textcolor{red}{\bf To be modified}

We suppose that $u^{\alpha,n}$ is bounded in space by $m^{\alpha}>-\infty$ and
$M^{\alpha}<+\infty$, and we denote 
$$\mathcal{U} = \prod_{\alpha=1}^{d}{[m^{\alpha},M^{\alpha}]} \quad \mbox{and} \quad \Lambda^{\alpha} = \sup_{\vect{u}\in\mathcal{U}}{\left|\lambda^\alpha(\vect{u})\right|}.$$ 
We say that $\vect{u}^{n}\in\mathcal{U}^{\Z}$ if $\vect{u}^{n}_i\in\mathcal{U}$ for all $i\in \Z$.
We now introduce two CFL conditions:
 \begin{equation}
    \frac{\Delta t}{\Delta x} \Lip(\vect{\lambda})TV(\vect{u}^n)<1
    \label{eqn:CFLexistence}
  \end{equation}
  \begin{equation}
    \frac{\Delta t}{\Delta x} \sum_{\alpha=1}^{d}{\Lambda^\alpha}<\frac{1}{2}
    \label{eqn:CFL}
  \end{equation}
We first prove that the semi-explicit scheme has a unique bounded solution at each time-step.
\begin{theorem}{\bf (Resolution of the semi-explicit scheme on one time step)}\label{thrm::1}
Assume that $\vect{\lambda}$ satisfies (\ref{eq::r1}). Let  $\vect{u}^{n}\in\mathcal{U}^{\Z}$, and assume that the two CFL conditions (\ref{eqn:CFLexistence}) and (\ref{eqn:CFL}) are satisfied.\\
 {\bf i) (Existence)}\\
Then there exists a unique solution $\vect{u}^{n+1}\in\mathcal{U}^{\Z}$ to the
  semi-explicit scheme (\ref{eqn:schema}).\\
  {\bf ii) (Monotonicity)}\\
Moreover if $\vect{u}^n$ is non-decreasing, i.e. satisfies
\begin{equation*}\label{eq::r13}
u^{n,\alpha}_{i+1}\ge u^{n,\alpha}_{i} \quad \mbox{for all}\quad i\in \Z,\quad \alpha = 1,...,d.
\end{equation*}
then $\vect{u}^{n+1}$ is also non-decreasing.
 \end{theorem}

\begin{remark}
  The resolution of the nonlinear problem boils down to the resolution
  of a local fixed point problem at each point $x_i=i\Delta x$. Note also that condition (\ref{eqn:CFLexistence}) is satisfied for $\vect{u}^{n}\in\mathcal{U}^{\Z}$ non decreasing if we have
\begin{equation*}\label{eq::r8}
 \frac{\Delta t}{\Delta x} \Lip(\vect{\lambda})\left(\sum_{\alpha=1,...,d}|M^\alpha-m^\alpha|\right)<1.
\end{equation*}
\end{remark}

Denoting $f(x)=x\ln(x)$, we then prove the following gradient entropy decay:
\begin{theorem}{\bf (Gradient entropy decay)}
  \label{thm:convexite}
  Assume that $\vect{\lambda}$ satisfies assumptions (\ref{eq::r1}) and (\ref{eq::7}).
Let us consider an initial data $\vect{u}^{0}\in\mathcal{U}$ which is assumed to be non-decreasing, i.e.
\begin{equation*}\label{eq::r7}
u^{0,\alpha}_{i+1}\ge u^{0,\alpha}_{i} \quad \mbox{for all}\quad i\in \Z,\quad \alpha = 1,...,d.
\end{equation*}
and let us consider the solution $\vect{u}^n$ of scheme (\ref{eqn:schema}), assuming the CFL conditions (\ref{eqn:CFLexistence}) and (\ref{eqn:CFL}) for all $n\ge 0$.
  Then $\theta_{i+\frac{1}{2}}^{\alpha,n}$ (defined in (\ref{eq::r9})) is non-negative for all
  $n\in\mathbb{N}$, and satisfies the following gradient entropy inequality for all $i_0\in \Z$ and $N\in \N$:
  \begin{equation}\label{eq:entropie}
    \sum_{\alpha=1}^{d}{\sum_{i\in I_N(i_0)}{f(\theta_{i+\frac{1}{2}}^{\alpha,n+1})}} \leq
    \sum_{\alpha=1}^{d}{\sum_{i\in I_N(i_0)}{f(\theta_{i+\frac{1}{2}}^{\alpha,n})}} 
    -\frac{\Delta t}{\Delta x}\sum_{\alpha=1}^{d}\left(F^{\alpha,n}_{i_0+N+1}-F^{\alpha,n}_{i_0-N}\right)
  \end{equation}
  where $I_N(i_0)$ is defined in (\ref{eq::r11}) and with the entropy flux
  \begin{equation}\label{eq:fluxentropie}
    F^{\alpha,n}_{i}=(\lambda^{\alpha,n+1}_i)_+f(\theta^{\alpha,n}_{i-\frac12})
    -(\lambda^{\alpha,n+1}_i)_-f(\theta^{\alpha,n}_{i+\frac12})
  \end{equation}
\end{theorem}
In particular, formally for $N=+\infty$, the flux terms on the boundary disappear on the right hand side of (\ref{eq:entropie}),
and we recover (\ref{eq::r10}).

We remark that $f$ can become negative. In order to ensure that every
term of the sum is nonnegative, we define $\tilde{f}$ as follows:
\begin{equation}\label{eq:ftilde}
  \forall \theta\geq 0,\quad \tilde{f}(\theta) = \left(f(\theta)+\frac{1}{e}\right)\char_{\{\theta>\frac{1}{e}\}}(\theta)
\end{equation}
$\tilde{f}$ is continuous, convex and nonnegative. A technical entropy
estimate similar to (\ref{eq:entropie})
is obtained on $\tilde{f}$ in Proposition \ref{prop:LlogL}, and will be
used to estimate the $L\log L$ norm of $\vect{u}_x$ at a discrete level.

%For a given $(n_0,i_0)\in\N\times \Z$, we will consider solutions $v$ to the linear scheme:
%\begin{equation}\label{eq::10}
%\frac{v^{n+1,\alpha}_{i}-v^{n,\alpha}_{i}}{\Delta t} 
%+ \left(\lambda^\alpha(u^{n_0}_{i_0})\right)_+ \left(\frac{v^{n,\alpha}_{i+1}-v^{n,\alpha}_{i}}{\Delta x}\right)
%+ \left(\lambda^\alpha(u^{n_0}_{i_0})\right)_- \left(\frac{v^{n,\alpha}_{i}-v^{n,\alpha}_{i-1}}{\Delta x}\right)=0,\quad \quad \mbox{for}\quad i\in\Z,\quad  n\ge n_0
%\end{equation}
%with data for $n=n_0$
%\begin{equation}\label{eq::11}
%v^{n_0}=u^{n_0}
%\end{equation}

% \begin{proposition}{\bf (A tame estimate for the scheme)}\label{pro::3}\\
% Under assumptions ({\bf TO PRECISE}), the following holds.
% There exists a constant $C>0$ ({\bf TO PRECISE}) such that for any $(n_0,i_0)\in\N\times \Z$ 
% and any integers $0< k \le  N$, we have:
% \begin{equation}\label{eq::9}
% \frac{1}{k\Delta t} \sum_{\alpha=1,...,d}\sum_{i\in I_{N-k}(i_0)}|u^{n_0+k,\alpha}_i-v^{n_0+k,\alpha}_i| 
% \Delta x \le C \left(TV[u^{n_0};I_{N}(i_0)]\right)^2
% \end{equation}
% where $v^n$ is solution to the scheme (\ref{eq::10}) with data for $n=n_0$ given in (\ref{eq::11}).
% \end{proposition}

Then we have the following
\begin{theorem}{\bf (Convergence of the solution of the scheme)}\label{th::2}\\
Assume that initial data $\vect{u}_0$ satisfies (\ref{eq::r2}), and that $\vect{\lambda}$ satisfies (\ref{eq::r1}) and (\ref{eq::7}).\\
Then there exists a bounded set $\mathcal{U}$ such that $\vect{u}_0(x)\in \mathcal{U}$ for all $x\in \R$.
Let us set the initial condition for the scheme
$$\vect{u}^0_i=\vect{u}_0(i\Delta x)$$
and let us consider the solution $(\vect{u}^n)_{n\ge 0}$ of the scheme (\ref{eqn:schema})
for time step $\Delta t >0$ and space step $\Delta x>0$ such that the CFL conditions (\ref{eqn:CFLexistence}) and (\ref{eqn:CFL})
are satisfied for all $n\ge 0$.
Let us call $\varepsilon=(\Delta t, \Delta x)$ and $\vect{u}^\varepsilon$ the function defined by
$$\vect{u}^\varepsilon(n \Delta t, i \Delta x) = \vect{u}^{n}_i \quad \mbox{for}\quad n\in\N,\quad i\in\Z$$
Then as $\varepsilon$ goes to zero, we have the following.\\
\noindent {\bf i) (Convergence for a subsequence)}\\ 
Up to extraction of a subsequence, there exists a continuous vanishing
viscosity solution $u$ of (\ref{eqn:system})-(\ref{eq::2}), such that
for any compact $K\subset [0,+\infty)\times \R$, we have 
$$|\vect{u}^\varepsilon-\vect{u}|_{L^\infty(K\cap (\Delta t \N)\times (\Delta x \Z),\R^d)} \to 0 \quad \mbox{as}\quad \varepsilon \to (0,0)$$
\noindent {\bf ii) (Convergence of the whole sequence)}\\ 
If we assume moreover that $\vect{\lambda}$ satisfies the strict hyperbolicity condition (\ref{eq::8}), then
the whole sequence $\vect{u}^\varepsilon$ converges to the unique continuous vanishing viscosity solution $\vect{u}$ of (\ref{eqn:system})-(\ref{eq::2}), as $\varepsilon$ goes to zero.
\end{theorem}

\begin{remark}
It would be interesting to adapt and extend the theory to the case where $\vect{\lambda}$ also depends on $(t,x)$.
At least for the scheme, this is an easy adaptation to write it.
It would also be interesting to extend the convergence of the solution of the scheme under the assumption of strict hyperbolicity  (\ref{eq::8}) without assuming (\ref{eq::7}) as a discrete analogue of Theorem 1.1 in \cite{EM2}.
\end{remark}

\begin{remark}
Indeed, we show a slightly better estimate than (\ref{eq:tameest}) without the ``limsup'', with explicit constants.
\end{remark}

\subsection{Literature}

For references on hyperbolic systems in non-conservative form, we refer to the references cited in \cite{EM1,EM2}.
Numerical schemes for hyperbolic systems are mainly written for systems in conservative form which enable to recover the correct Rankine-Hugoniot shock relations. We refer to \cite{leveque} for a review of the main classes of existing schemes. Among these schemes, convergence results are seldom found for hyperbolic systems. 

The Lax-Wendroff theorem \cite{LaxWendroff} shows that if a consistent and conservative numerical scheme converges (in $L^1$ with bounded total variation), its limit is a weak solution to the hyperbolic system. However, in order to obtain convergence of the scheme, stability is needed, in general in the form of TV-stability. For the scalar Godunov scheme, convergence is obtained due to its total variation diminishing (TVD) property. This is no longer the case for systems \cite{leveque}. Stability can still be proved for certain special systems of two equations, for instance in \cite{temple1,temple2,levequetemple}. Similar results can be obtained for a class of nonlinear systems with straight-line fields \cite[pp.~102--103]{bressan}. Nonlinear stability can also be assessed through the use of invariant domains and entropy inequalities \cite{bouchut}, for HLL, HLLC and kinetic solvers for Euler equations of gas dynamics.

In the case of conservative systems where the initial data has sufficiently small total variation, Glimm's random choice method \cite{glimm} is provably convergent. A deterministic variant (replacing random with equidistributed sampling) has also been proved to converge under the same assumptions \cite{liu}. We are not aware of convergence results of numerical schemes for non-conservative hyperbolic systems with large initial data.

\subsection{Outline of the article}

This paper is organized as follows. In Section \ref{sec:preliminary}, we prove some preliminary results on the existence of a solution to the scheme (Theorem \ref{thrm::1}), on the monotonicity and boundedness of the solution, and a discrete analogue of the tame estimate given in Definition \ref{defi::1}. We then prove the decrease of the discrete entropy (Theorem \ref{thm:convexite}) in Section \ref{sec:entropy}. In addition, we establish a similar entropic estimate for the scheme. Finally, in Section \ref{sec:convergence}, we sum up all the results and prove the convergence of the scheme (Theorem \ref{th::2}).

\section{Preliminary results on the scheme}
\label{sec:preliminary}

\subsection{Existence and uniqueness of the solution of the semi-explicit
  scheme}

\begin{Proofc}{Proof of Theorem \ref{thrm::1}, part i)}
  We define the truncature $T\lambda^\alpha$ of $\lambda^\alpha$ by $\Lambda^\alpha$:
  \begin{equation*}
    T\lambda^\alpha(\vect{u}) = \left\{
      \begin{array}{cl}
        \lambda^\alpha(\vect{u}) & \text{if }
        |\lambda^\alpha(\vect{u})|\leq\Lambda^\alpha \\
        \Lambda^{\alpha} & \text{if }
        \lambda^\alpha(\vect{u})>\Lambda^\alpha \\
        -\Lambda^{\alpha} & \text{if }
        \lambda^\alpha(\vect{u})<-\Lambda^\alpha
      \end{array}\right.
  \end{equation*}
  $\vect{T\lambda}$ is also Lipschitz and $\Lip(\vect{T\lambda})\leq\Lip(\vect{\lambda})$.
  For $\vect{v}\in \R^d$, let us define the function
  $\vect{F}_{\vect{u}_i^n,\vect{u}_{i-1}^n,\vect{u}_{i+1}^n}$ such that, for all $\alpha\in\{1,\dots,d\}$,
  \begin{equation*}
    F^\alpha_{\vect{u}_i^n,\vect{u}_{i-1}^n,\vect{u}_{i+1}^n}(\vect{v})
    = u_{i}^{\alpha,n}+\frac{\Delta t}{\Delta x}\left((T\lambda^\alpha(\vect{v}))_-(u_{i+1}^{\alpha,n}-u_{i}^{\alpha,n})-(T\lambda^\alpha(\vect{v}))_+(u_{i}^{\alpha,n}-u_{i-1}^{\alpha,n})\right) 
  \end{equation*}
  Then the scheme  (\ref{eqn:schema}) can be written (if $\vect{u}^n\in \mathcal{U}^{\Z}$) as
  \begin{equation}\label{eq::r12}
\vect{u}_i^{n+1}=  \vect{F}_{\vect{u}_i^n,\vect{u}_{i-1}^n,\vect{u}_{i+1}^n}(\vect{u}^{n+1}_i)
\end{equation}
  We observe that,
  for all $\vect{u}$ and $\vect{v}$ in $\R^d$, for all $\alpha\in\{1,\dots,d\}$,
  \begin{align*}
    |
    F^\alpha_{\vect{u}_i^n,\vect{u}_{i-1}^n,\vect{u}_{i+1}^n}(\vect{u})-
    F^\alpha_{\vect{u}_i^n,\vect{u}_{i-1}^n,\vect{u}_{i+1}^n}(\vect{v})|
    & \leq\frac{\Delta
      t}{\Delta
      x}|T\lambda^\alpha(\vect{u})-T\lambda^\alpha(\vect{v})|\left(|u_{i+1}^{\alpha,n}-u_{i}^{\alpha,n}|+|u_{i}^{\alpha,n}-u_{i-1}^{\alpha,n}|\right)\\
    & \leq \frac{\Delta t}{\Delta x}\Lip(\vect{\lambda})TV(\vect{u}^{n})|\vect{u}-\vect{v}|
  \end{align*}
  So that $\vect{F}_{\vect{u}_i^n,\vect{u}_{i-1}^n,\vect{u}_{i+1}^n}$ is
  contractive on $\R^d$ thanks to CFL condition (\ref{eqn:CFLexistence}) and Banach fixed point theorem
  yields the existence and uniqueness of a solution $\vect{u}$ of (\ref{eq::r12}) on $\R^d$. 
  
  In addition, due to the CFL condition (\ref{eqn:CFL}),
  $F^\alpha_{\vect{u}_i^n,\vect{u}_{i-1}^n,\vect{u}_{i+1}^n}(\vect{u})$
  is a convex combination of the $\vect{u}_i^n$, $\vect{u}_{i-1}^n$ and
  $\vect{u}_{i+1}^n$ contained in the convex $\mathcal{U}$, so that
  $\vect{u} =
  F^\alpha_{\vect{u}_i^n,\vect{u}_{i-1}^n,\vect{u}_{i+1}^n}(\vect{u})$
  is also in $\mathcal{U}$. As $\vect{T\lambda}=\vect{\lambda}$ on
  $\mathcal{U}$, we can conclude that the unique fixed point  of equation (\ref{eq::r12})
  is the solution $\vect{u}_i^{n+1}$ of the scheme (\ref{eqn:schema}).
\end{Proofc}

\subsection{Expression of $\theta^{n+1}$}

We derive an equation for the evolution of
$\theta_{i+\frac{1}{2}}^{\alpha,n}$ in time.

\begin{lemma}{\bf (Evolution of $\theta$)}
  \label{lem:theta}
  Let $u_i^{\alpha,n}$ be the solution of the semi-explicit scheme
  (\ref{eqn:schema}).
  Then $\theta_{i+\frac{1}{2}}^{\alpha,n+1}$ satisfies the following
  relation:
  \begin{multline}
    \theta_{i+\frac{1}{2}}^{\alpha,n+1} =
    \left(1-\frac{\Delta t}{\Delta
        x}((\lambda_{i+1}^{\alpha,n+1})_++(\lambda_{i}^{\alpha,n+1})_-)\right)\theta_{i+\frac{1}{2}}^{\alpha,n}\\
    + \frac{\Delta t}{\Delta
      x}(\lambda_{i+1}^{\alpha,n+1})_-\theta_{i+\frac{3}{2}}^{\alpha,n}
    + \frac{\Delta t}{\Delta
      x}(\lambda_{i}^{\alpha,n+1})_+\theta_{i-\frac{1}{2}}^{\alpha,n}
    \label{eqn:theta}
  \end{multline}
\end{lemma}

\begin{Proof}
  With the definition of $\theta_{i+\frac{1}{2}}^{\alpha,n}$, we observe:
  \begin{equation*}
    \theta_{i+\frac{1}{2}}^{\alpha,n+1} = \theta_{i+\frac{1}{2}}^{\alpha,n} +
    \frac{\Delta t}{\Delta x}\frac{u_{i+1}^{\alpha,n+1}-u_{i+1}^{\alpha,n}}{\Delta t} -
    \frac{\Delta t}{\Delta
      x}\frac{u_{i}^{\alpha,n+1}-u_{i}^{\alpha,n}}{\Delta
      t} 
  \end{equation*}
  Inserting (\ref{eqn:schema}) at points $x_i$ and $x_{i+1}$, we get
  equation (\ref{eqn:theta}).
\end{Proof}

\subsection{$\vect{u}^{n}$ is non-decreasing if
  $\vect{u}^{0}$ is non-decreasing}

\begin{lemma}{\bf (Monotonicity)}
  \label{lem:croissance}
  Let $\vect{u}^{n}\in\mathcal{U}^{\Z}$ be non-decreasing. 
  Assume that $\vect{\lambda}$ satisfies (\ref{eq::r1}) and assume the two CFL conditions 
  (\ref{eqn:CFLexistence}) and (\ref{eqn:CFL}).
   Then $\vect{u}^{n+1}$ is non-decreasing.
\end{lemma}

\begin{Proof}
  In equation (\ref{eqn:theta}), the coefficients $\frac{\Delta t}{\Delta
    x}(\lambda_{i+1}^{\alpha,n+1})_-$ and $\frac{\Delta t}{\Delta
    x}(\lambda_{i}^{\alpha,n+1})_+$ are positive by definition,
  Theorem \ref{thrm::1}, part i) yields that $\vect{u}^{n+1}$ is in
  $\mathcal{U}$ and using the CFL condition (\ref{eqn:CFL}), we obtain that:
  \begin{equation*}
    \left(1-\frac{\Delta t}{\Delta
        x}((\lambda_{i+1}^{\alpha,n+1})_++(\lambda_{i}^{\alpha,n+1})_-)\right)\geq 0
  \end{equation*}
  As $u_i^{\alpha,n}$ is non-decreasing, for all $i\in\mathbb{N}$ and
  $1\leq\alpha\leq d$,
  $\theta_{i+\frac{1}{2}}^{\alpha,n}\geq 0$, and therefore
  $\theta_{i+\frac{1}{2}}^{\alpha,n+1}$ is non-negative too. This is
  equivalent to $u_i^{\alpha,n+1}$ non-decreasing.
\end{Proof}

\begin{Proofc}{Proof of Theorem \ref{thrm::1}, part ii)}
We simply apply a recursion on $n\ge 0$, using Lemma \ref{lem:croissance}.
\end{Proofc}

\subsection{$\vect{u}$ has a non-increasing total variation}

\begin{lemma}{\bf (Total variation decay)}
  \label{lem:TVD}
  Let $\vect{u}^{n}\in\mathcal{U}^{\Z}$.
   Assume that $\vect{\lambda}$ satisfies (\ref{eq::r1}) and assume the two CFL conditions 
  (\ref{eqn:CFLexistence}) and (\ref{eqn:CFL}).
   Let $i_0\in \Z$ be a fixed index and $N\in\mathbb{N}\backslash \{0\}$. Then:
   \begin{equation*}
     TV[\vect{u}^{n+1};I_{N-1}(i_0)]\leq TV[\vect{u}^{n};I_N(i_0)]
   \end{equation*}
   and
   \begin{equation}\label{eq:TVdecay}
     TV(\vect{u}^{n+1})\leq TV(\vect{u}^{n}) \text{ if } TV(\vect{u}^{n})<+\infty
   \end{equation}
\end{lemma}
\begin{Proof}
  CFL condition (\ref{eqn:CFL}) allows us to write
  $u_i^{\alpha,n+1}$ as a convex sum of $u_{i-1}^{\alpha,n}$,
  $u_{i}^{\alpha,n}$ and $u_{i+1}^{\alpha,n}$, so that:
  \begin{multline*}
    |u_{i+1}^{\alpha,n+1}-u_i^{\alpha,n+1}|\leq\left(1-\frac{\Delta
        t}{\Delta
        x}[(\lambda_i^{\alpha,n+1})_-+(\lambda_{i+1}^{\alpha,n+1})_+]\right)|u_{i+1}^{\alpha,n}-u_i^{\alpha,n}|  \\
    + \frac{\Delta t}{\Delta
      x}(\lambda_{i+1}^{\alpha,n+1})_-|u_{i+2}^{\alpha,n}-u_{i+1}^{\alpha,n}|
    + \frac{\Delta t}{\Delta
      x}(\lambda_{i}^{\alpha,n+1})_+|u_{i}^{\alpha,n}-u_{i-1}^{\alpha,n}|
  \end{multline*}
  Summing these terms for $i\in I_{N-1}(i_0)$ gives a sum for $i\in I_{N-1}(i_0)$ of
  $|u_{i+1}^{\alpha,n}-u_i^{\alpha,n}|$, and the remaining terms are for 
  $i\in I_N(i_0)\backslash  I_{N-1}(i_0)$ with coefficients inferior to 1 due to CFL condition (\ref{eqn:CFL}). 
\end{Proof}

\subsection{A tame estimate for the scheme}

In this subsection, we prove a discrete analogue to the continuous
vanishing viscosity solution given in Definition \ref{defi::1} for the
discrete solution $u_i^n$.

\begin{proposition}{\bf (Discrete tame estimate)}
 Let $\vect{u}^{n}\in\mathcal{U}^{\Z}$. Assume that $\vect{\lambda}$ satisfies (\ref{eq::r1}) and assume the two CFL conditions 
  (\ref{eqn:CFLexistence}) and (\ref{eqn:CFL}).
  Then the following holds.
  Let $(i_0,n_0)\ \Z\times \N$ be a fixed. Let $(\vect{v}^n)_{n\ge n_0}$ be
  the solution of the explicit discretization of the linear hyperbolic
  Cauchy problem with frozen constant coefficients for $n\ge n_0$:
  \begin{equation}\label{eq::r15}
    \frac{v_i^{\alpha,n+1}-v_i^{\alpha,n}}{\Delta t}
    -(\lambda^{\alpha}(\vect{u}_{i_0}^{n_0}))_-\left(\frac{v_{i+1}^{\alpha,n}-v_i^{\alpha,n}}{\Delta
    x}\right)+(\lambda^{\alpha}(\vect{u}_{i_0}^{n_0}))_+\left(\frac{v_{i}^{\alpha,n}-v_{i-1}^{\alpha,n}}{\Delta
    x}\right) = 0
  \end{equation}
  with $\vect{v}^{n_0} = \vect{u}^{n_0}$. Then, for all $k\in\mathbb{N}\backslash \{0\}$ such
  that $k\le N$, 
  \begin{equation}\label{eq::9}
    \frac{1}{k\Delta
      t}\sum_{\alpha=1}^d\sum_{i\in
      I_{N-k}(i_0)}{|u_i^{\alpha,n_0+k}-v_i^{\alpha,n_0+k}|\Delta x}\leq
    2\Lip(\vect{\lambda})\left(TV[u^{n_0};I_N(i_0)]\right)^2.
  \end{equation}
\end{proposition}
\begin{Proof}
  Let:
  \begin{equation*}
    \mathcal{I}_N^{k} = \sum_{\alpha=1}^d{\sum_{i\in I_N(i_0)}{|u_i^{\alpha,n_0+k}-v_i^{\alpha,n_0+k}|}}
  \end{equation*}
  Using the schemes (\ref{eqn:schema}) and (\ref{eq::r15}), we obtain:
  \begin{align*}
    \mathcal{I}_{N-k-1}^{k+1}&\leq\sum_{\alpha=1}^d{\sum_{i\in I_{N-k-1}(i_0)}{\left|(1-\frac{\Delta
        t}{\Delta x}|\lambda^\alpha(\vect{u}_{i_0}^{n_0})|)(u_i^{\alpha,n_0+k}-v_i^{\alpha,n_0+k})\right|}}\\
        &+\sum_{\alpha=1}^d{\sum_{i\in I_{N-k-1}(i_0)}{\frac{\Delta t}{\Delta x}(\lambda^\alpha(\vect{u}_{i_0}^{n_0}))_-\left|u_{i+1}^{\alpha,n_0+k}-v_{i+1}^{\alpha,n_0+k}\right|}}\\
        &+\sum_{\alpha=1}^d{\sum_{i\in I_{N-k-1}(i_0)}{\frac{\Delta t}{\Delta x}(\lambda^\alpha(\vect{u}_{i_0}^{n_0}))_+\left|u_{i-1}^{\alpha,n_0+k}-v_{i-1}^{\alpha,n_0+k}\right|}}\\
        &+\sum_{\alpha=1}^d{\sum_{i\in I_{N-k-1}(i_0)}{\left|\frac{\Delta t}{\Delta x}((\lambda^\alpha(\vect{u}_{i}^{n_0+k+1}))_--(\lambda^\alpha(\vect{u}_{i_0}^{n_0}))_-)(u_{i+1}^{\alpha,n_0+k}-u_{i}^{\alpha,n_0+k})\right|}}\\
        &+\sum_{\alpha=1}^d{\sum_{i\in I_{N-k-1}(i_0)}{\left|\frac{\Delta t}{\Delta x}((\lambda^\alpha(\vect{u}_{i}^{n_0+k+1}))_+-(\lambda^\alpha(\vect{u}_{i_0}^{n_0}))_+)(u_{i}^{\alpha,n_0+k}-u_{i-1}^{\alpha,n_0+k})\right|}}
  \end{align*}
  CFL condition (\ref{eqn:CFL}) gives us that $1-\frac{\Delta t}{\Delta x}|\lambda^\alpha(\vect{u}_{i_0}^{n_0})|$ is positive. 
  In the right hand side of the inequality, the first three terms can then be controlled by $$\sum_{\alpha=1}^d{\sum_{i\in I_{N-k}(i_0)}{\left|u_i^{\alpha,n_0+k}-v_i^{\alpha,n_0+k}\right|}}=\mathcal{I}_{N-k}^{k}$$.

  We note that, for all $\vect{u}, \vect{v}\in\mathbb{R}^d$,
    $|(\lambda^\alpha(\vect{u}))_--(\lambda^\alpha(\vect{v}))_-|\leq|\lambda^\alpha(\vect{u})-\lambda^\alpha(\vect{v})|$
    and
    $|(\lambda^\alpha(\vect{u}))_+-(\lambda^\alpha(\vect{v}))_+|\leq|\lambda^\alpha(\vect{u})-\lambda^\alpha(\vect{v})|$,
    and we recall that: 
	\begin{equation*}
		|\lambda^\alpha(\vect{u}_{i}^{n_0+k+1})-\lambda^\alpha(\vect{u}_{i_0}^{n_0})|\leq\Lip(\vect{\lambda})|\vect{u}_{i}^{n_0+k+1}-\vect{u}_{i_0}^{n_0}|
	\end{equation*}
    Using the same convexity argument as in Lemma
    \ref{lem:croissance}, it is easy to see that if, for some
    $K\in\mathbb{N}\backslash \{0\}$, we have
    $$m_n^{\alpha}(I_K(i_0))\le u^{\alpha,n} \le M_n^{\alpha}(I_K(i_0)),\quad\mbox{for all}\quad i\in I_K(i_0)$$
    then we have
    $$m_n^{\alpha}(I_K(i_0))\le u^{\alpha,n+1} \le M_n^{\alpha}(I_K(i_0)),\quad\mbox{for all}\quad i\in I_{K-1}(i_0)$$
    A straightforward recursion yields that
    $u^{\alpha,n_0+k+1}$ is bounded on $I_{N-k-1}(i_0)$ by the bounds of
    $u^{\alpha,n_0}$ on $I_N(i_0)$,
    $m_{n_0}^{\alpha}(I_N(i_0))$ and $M_{n_0}^{\alpha}(I_N(i_0))$. As a result,
    for all $i\in I_{N-k-1}(i_0)$, 
    \begin{align*}
      |\lambda^\alpha(\vect{u}_{i}^{n_0+k+1})-\lambda^\alpha(\vect{u}_{i_0}^{n_0})|&\leq\Lip(\vect{\lambda})|\vect{M}_{n_0}(I_N(i_0))-\vect{m}_{n_0}(I_N(i_0))|\\
      &\leq\Lip(\vect{\lambda})TV[\vect{u}^{n_0};I_N(i_0)]
    \end{align*}
In the end, using Lemma \ref{lem:TVD}, we deduce that
\begin{equation*}
  \mathcal{I}_{N-k-1}^{k+1}\leq \mathcal{I}_{N-k}^{k}+2\Delta t\Lip(\vect{\lambda})\left(TV[\vect{u}^{n_0};I_N(i_0)]\right)^2
\end{equation*}
The result is then obtained through a straightforward recursion on $k$.
    
\end{Proof}

%%%%%%%%%%%%%%%%%%%%%%%%%%
\section{The gradient entropy}
\label{sec:entropy}

In this section, we define $f:\mathbb{R}^+\rightarrow\mathbb{R}$ the
convex function $f(x) = x\ln(x)$. 

\subsection{Preparatory lemma}

\begin{lemma}{\bf (Convexity inequality for $f$)}
  \label{lem:convexity}
  Let $a_k$ and $\theta_k$ be two finite sequences of non-negative real numbers such that
  $0<\sum_{k}{a_k}<+\infty$. Define:
  \begin{equation*}
    \theta = \sum_{k}{a_k\theta_{k}}
  \end{equation*}
  Then the following inequality holds:
  \begin{equation*}
    f(\theta) \leq \sum_{k}{a_kf(\theta_{k})} + \theta\ln\left(\sum_{k}{a_k}\right)
  \end{equation*}
\end{lemma}

\begin{Proof}
  As $\sum_k{a_k}>0$,
  \begin{equation*}
    \frac{1}{\sum_k{a_k}}\theta = \sum_k{\frac{a_k}{\sum_l{a_l}}\theta_{k}}
  \end{equation*}
  is a convex sum of the $\theta_{k}\geq 0$. Using the convexity of
  $f$ on $\mathbb{R}^+$, 
  \begin{equation*}
    f\left(\frac{1}{\sum_k{a_k}}\theta\right) \leq \sum_k{\frac{a_k}{\sum_l{a_l}}f(\theta_{k})}
  \end{equation*}
  Using the expression of $f(x)=x\ln(x)$, 
  \begin{equation*}
    f\left(\frac{1}{\sum_k{a_k}}\theta\right) = \frac{1}{\sum_k{a_k}}\left(f(\theta) - \theta\ln\left(\sum_k{a_k}\right)\right)
  \end{equation*}
  which proves the result.
\end{Proof}

\subsection{Proof of Theorem \ref{thm:convexite}}

\begin{Proof}
  For $i$, $\alpha$ and $n$ fixed, Lemma \ref{lem:theta} gives us the expression
  (\ref{eqn:theta}) for $\theta_{i+\frac{1}{2}}^{\alpha,n+1}$. Let us
  remark that the coefficients $a_1=\frac{\Delta t}{\Delta
    x}(\lambda_{i+1}^{\alpha,n+1})_-$ and $a_2=\frac{\Delta t}{\Delta
    x}(\lambda_{i}^{\alpha,n+1})_+$ are non-negative by definition, and that
  CFL condition (\ref{eqn:CFL}) yields:
  \begin{equation*}
    a_3=\left(1-\frac{\Delta t}{\Delta
        x}((\lambda_{i+1}^{\alpha,n+1})_++(\lambda_{i}^{\alpha,n+1})_-)\right)
    \geq 0
  \end{equation*}
  Defining $\mu_{i+\frac{1}{2}}^{\alpha,n+1}=1-(a_1+a_2+a_3)$, let us note that CFL condition (\ref{eqn:CFLexistence}) joined to (\ref{eq:TVdecay}) also gives:
  \begin{equation}\label{eq:r48}
    1-\mu_{i+\frac{1}{2}}^{\alpha,n+1}=a_1+a_2+a_3=\left(1-\frac{\Delta t}{\Delta
        x}(\lambda_{i+1}^{\alpha,n+1}-\lambda_{i}^{\alpha,n+1})\right) >0
  \end{equation}
  Using Lemma \ref{lem:convexity} on the convex sum, we obtain:
  \begin{multline*}
    f(\theta_{i+\frac{1}{2}}^{\alpha,n+1}) \leq
    \left(1-\frac{\Delta t}{\Delta
        x}((\lambda_{i+1}^{\alpha,n+1})_++(\lambda_{i}^{\alpha,n+1})_-)\right)f(\theta_{i+\frac{1}{2}}^{\alpha,n})\\
    + \frac{\Delta t}{\Delta
      x}(\lambda_{i+1}^{\alpha,n+1})_-f(\theta_{i+\frac{3}{2}}^{\alpha,n})
    + \frac{\Delta t}{\Delta
      x}(\lambda_{i}^{\alpha,n+1})_+f(\theta_{i-\frac{1}{2}}^{\alpha,n})\\
    + \theta_{i+\frac{1}{2}}^{\alpha,n+1}\ln(1-\mu_{i+\frac{1}{2}}^{\alpha,n+1})
  \end{multline*}
  Arranging terms, the expression exhibits a discrete divergence form:
  \begin{multline}
    f(\theta_{i+\frac{1}{2}}^{\alpha,n+1}) \leq
    f(\theta_{i+\frac{1}{2}}^{\alpha,n})
    + \frac{\Delta t}{\Delta
      x}\left((\lambda_{i+1}^{\alpha,n+1})_-f(\theta_{i+\frac{3}{2}}^{\alpha,n})-(\lambda_{i}^{\alpha,n+1})_-f(\theta_{i+\frac{1}{2}}^{\alpha,n})\right)\\
    -\frac{\Delta t}{\Delta
      x} \left((\lambda_{i+1}^{\alpha,n+1})_+f(\theta_{i+\frac{1}{2}}^{\alpha,n})-((\lambda_{i}^{\alpha,n+1})_+f(\theta_{i-\frac{1}{2}}^{\alpha,n})\right)\\
    +
    \theta_{i+\frac{1}{2}}^{\alpha,n+1}\ln(1-\mu_{i+\frac{1}{2}}^{\alpha,n+1})
    \label{eqn:divergence}
  \end{multline}
  Summing (\ref{eqn:divergence}) over $i\in I_N(i_0)$ and over $\alpha$,
  the second and third terms cancel and we obtain:
  \begin{multline}
    \sum_{\alpha=1}^{d}{\sum_{i\in I_N(i_0)}{f(\theta_{i+\frac{1}{2}}^{\alpha,n+1})}} \leq
    \sum_{\alpha=1}^{d}{\sum_{i\in I_N(i_0)}{f(\theta_{i+\frac{1}{2}}^{\alpha,n})}}
    +
    \sum_{\alpha=1}^{d}{\sum_{i\in I_N(i_0)}{\theta_{i+\frac{1}{2}}^{\alpha,n+1}\ln(1-\mu_{i+\frac{1}{2}}^{\alpha,n+1})}} \\
    - \frac{\Delta t}{\Delta x}\sum_{\alpha=1}^{d}{(F_{i_0+N+1}^{\alpha,n}-F_{i_0-N}^{\alpha,n})}
    \label{eq:ineglog}
  \end{multline}
  with $F^{\alpha,n}_{i}$ defined in (\ref{eq:fluxentropie}).
  We observe that $\ln(1-\mu)\leq-\mu$ for all $\mu< 1$, we note that
  $\mu_{i+\frac{1}{2}}^{\alpha,n+1}<1$ due to (\ref{eq:r48}) and we
  recall that $\theta_{i+\frac{1}{2}}^{\alpha,n+1}$ is non-negative, so
  that
  \begin{multline*}
    \sum_{\alpha=1}^{d}{\sum_{i\in I_N(i_0)}{f(\theta_{i+\frac{1}{2}}^{\alpha,n+1})}} \leq
    \sum_{\alpha=1}^{d}{\sum_{i\in I_N(i_0)}{f(\theta_{i+\frac{1}{2}}^{\alpha,n})}}
    -
    \sum_{\alpha=1}^{d}{\sum_{i\in I_N(i_0)}{\theta_{i+\frac{1}{2}}^{\alpha,n+1}\mu_{i+\frac{1}{2}}^{\alpha,n+1}}} \\
    -  \frac{\Delta t}{\Delta x}\sum_{\alpha=1}^{d}{(F_{i_0+N+1}^{\alpha,n}-F_{i_0-N}^{\alpha,n})}
  \end{multline*}
  Now by definition:
  \begin{align*}
    \mu_{i+\frac{1}{2}}^{\alpha,n+1} &=
    \frac{\Delta t}{\Delta
      x}\left(\lambda^\alpha(\vect{u}_{i+1}^{n+1})-\lambda^\alpha(\vect{u}_{i}^{n+1})\right) \\
    &=  \frac{\Delta t}{\Delta
      x}\sum_{\beta=1}^{d}{\int_0^1{\frac{\partial\lambda^\alpha}{\partial
            u^\beta}\left(\vect{u}_i^{n+1}+\tau(\vect{u}_{i+1}^{n+1}-\vect{u}_i^{n+1})\right)\cdot(u_{i+1}^{\beta,n+1}-u_{i}^{\beta,n+1})d\tau}}
  \end{align*}
  
  So that, summing over $\alpha$ and using the definition of $\theta_{i+\frac{1}{2}}^{\alpha,n}$:
  \begin{equation}
    \sum_{\alpha=1}^{d}{\mu_{i+\frac{1}{2}}^{\alpha,n+1}\theta_{i+\frac{1}{2}}^{\alpha,n+1}}
    =  \Delta
    t\int_0^1{\left(\matr{\nabla}\vect{\lambda}(\vect{u}_i^{n+1}+\tau\Delta x\vect{\theta}_{i+\frac{1}{2}}^{n+1})\cdot\vect{\theta}_{i+\frac{1}{2}}^{n+1}\right)\cdot\vect{\theta}_{i+\frac{1}{2}}^{n+1}d\tau}
    \geq 0 \label{eq:muthetapositive}
  \end{equation}
  where we have used assumption (\ref{eq::7}). In the end, we obtain the gradient entropy decay:
  \begin{equation*}
    \sum_{\alpha=1}^{d}{\sum_{i\in I_N(i_0)}{f(\theta_{i+\frac{1}{2}}^{\alpha,n+1})}} \leq
    \sum_{\alpha=1}^{d}{\sum_{i\in I_N(i_0)}{f(\theta_{i+\frac{1}{2}}^{\alpha,n})}} 
    -  \frac{\Delta t}{\Delta x}\sum_{\alpha=1}^{d}{(F_{i_0+N+1}^{\alpha,n}-F_{i_0-N}^{\alpha,n})}
  \end{equation*}
\end{Proof}

\subsection{Gradient entropy estimate}

As $f$ is negative for $\theta\in(0,\frac{1}{e})$, we use the following
similar result on $\tilde{f}$ as defined in (\ref{eq:ftilde}) in order to have a discrete estimate on $\vect{u}_x$ in the $L\log L$ norm:
\begin{proposition}{\bf (Gradient entropy estimate for the scheme)}\label{prop:LlogL}
  Under the assumptions of Theorem \ref{thm:convexite}, we have
  \begin{equation*}
    \sum_{\alpha=1}^{d}{\sum_{i\in\Z}{\tilde{f}(\theta_{i+\frac{1}{2}}^{\alpha,n+1})\Delta
      x}} \leq
    \sum_{\alpha=1}^{d}{\sum_{i\in\Z}{\tilde{f}(\theta_{i+\frac{1}{2}}^{\alpha,n})\Delta
      x}}
    + C\Delta t
  \end{equation*}
  if the right hand side is finite, with $C=C_2d\Lip(\vect{\lambda})TV(\vect{u}^0)$ where $C_2 = \frac{1}{e\ln 2}$.
\end{proposition}

In order to prove this result, we first need two technical lemmata on
$\tilde{f}$, analogous to Lemma \ref{lem:convexity}. 

\begin{lemma}{\bf (Technical estimate)}
  \label{lem:ftilde}
  Let $\gamma_m> 1$. There
  exists a non-negative function $g(\theta,\gamma)$ and a constant
  $C_{\gamma_m}>0$ (depending only on $\gamma_m$) such that, for all
  $\theta>0$ and $\gamma\in(0,\gamma_m)$, 
  \begin{equation}\label{eq:r58}
    \tilde{f}\left(\frac{\theta}{\gamma}\right)\geq \frac{1}{\gamma}\tilde{f}(\theta)-\frac{1}{\gamma}g(\theta,\gamma)\ln(\gamma)
  \end{equation}
  and:
  \begin{equation*}
    |\theta-g(\theta,\gamma)|\leq C_{\gamma_m} = \frac{\gamma_m-1}{e\ln(\gamma_m)}
  \end{equation*}
\end{lemma}
\begin{Proof}
  We detail the four cases:
  \begin{itemize}
  \item[\bf Case A:] $\frac{\theta}{\gamma}\geq\frac{1}{e}$ and $\theta\geq
    \frac{1}{e}$.\\
    We have for $\gamma\neq 1$:
    \begin{equation*}
      \frac{1}{\ln(\gamma)}\left(\tilde{f}(\theta)-\gamma\tilde{f}\left(\frac{\theta}{\gamma}\right)\right) = \theta - \frac{1}{e}\frac{\gamma-1}{\ln(\gamma)}
    \end{equation*}
    We then set for any $\gamma>0$:
    \begin{equation}\label{eq:r61}
      g(\theta,\gamma) = \left(\theta - \frac{1}{e}\frac{\gamma-1}{\ln(\gamma)}\right)_+
    \end{equation}
    This implies (\ref{eq:r58}) for $\gamma\geq 1$.
    As $\gamma\in(0,\gamma_m)$, and $\frac{\gamma-1}{\ln(\gamma)}$ is
    non-negative increasing,  we get
    \begin{equation*}
      |\theta-g(\theta,\gamma)|\leq \frac{1}{e}\frac{\gamma_m-1}{\ln(\gamma_m)}
    \end{equation*}
    Now for $\gamma\leq 1$, we have 
    $$g(\theta,\gamma)=\theta-\frac{1}{e}\frac{\gamma-1}{\ln(\gamma)}\geq g(\theta,1)=\theta-\frac{1}{e}\geq 0$$
    This shows that (\ref{eq:r58}) still holds for $0<\gamma\leq 1$.
  \item[\bf Case B:] $\frac{\theta}{\gamma}\geq\frac{1}{e}$ and $\theta<
    \frac{1}{e}$.\\
    Then we have $0<\gamma<1$ and
    \begin{align*}
      \tilde{f}\left(\frac{\theta}{\gamma}\right)-\frac{1}{\gamma}\tilde{f}(\theta)
      & =
      \frac{1}{\gamma}\theta\ln(\theta)-\frac{1}{\gamma}\theta\ln(\gamma)+\frac{1}{e}
      \\
      & \geq -\frac{1}{\gamma}\theta\ln(\gamma)+\frac{1}{e}\frac{\gamma-1}{\gamma} = -\frac{1}{\gamma}g(\theta,\gamma)\ln(\gamma)
    \end{align*}
    for $g(\theta,\gamma)$ defined in (\ref{eq:r61}).
  \item[\bf Case C:] $\frac{\theta}{\gamma}<\frac{1}{e}$ and $\theta\geq
    \frac{1}{e}$.
    \begin{equation*}
      \tilde{f}\left(\frac{\theta}{\gamma}\right)-\frac{1}{\gamma}\tilde{f}(\theta)
      \geq -\frac{1}{\gamma}\theta\ln(\gamma)
    \end{equation*}
    We take $g(\theta,\gamma) = \theta\geq 0$ in
    this case.
  \item[\bf Case D:] $\frac{\theta}{\gamma}<\frac{1}{e}$ and $\theta<
    \frac{1}{e}$. 
    \begin{equation*}
      \tilde{f}\left(\frac{\theta}{\gamma}\right)-\frac{1}{\gamma}\tilde{f}(\theta)
      \geq 0
    \end{equation*}
    We take $g(\theta,\gamma) = 0$ in this
    case, and we check that:
    \begin{equation*}
      |\theta-g(\theta,\gamma)|= \theta\leq \frac{1}{e}\leq \frac{\gamma_m-1}{e\ln(\gamma_m)}
    \end{equation*}
  \end{itemize}
\end{Proof}

\begin{lemma}{\bf (Convexity inequality for $\tilde{f}$)}
  \label{lem:convexitytilde}
  Let $a_k$ and $\theta_k$ be two finite sequences of non-negative real numbers such that
  $0<\sum_{k}{a_k}<2$. Define:
  \begin{equation*}
    \theta = \sum_{k}{a_k\theta_{k}}
  \end{equation*}
  Then the following inequality holds:
  \begin{equation*}
    \tilde{f}(\theta) \leq \sum_{k}{a_k\tilde{f}(\theta_{k})} + g\left(\theta,\sum_{k}{a_k}\right)\ln\left(\sum_{k}{a_k}\right)
  \end{equation*}
  where $g(\theta,\gamma)$ is given by Lemma
  \ref{lem:ftilde} for $\gamma_m=2$. 
\end{lemma}
\begin{Proof}
  The proof is an adaptation of the proof of Lemma \ref{lem:convexity} with $\gamma = \sum_k{a_k}$,
  using Lemma   \ref{lem:ftilde} in the convexity inequality for $\tilde{f}$. 
\end{Proof}

\begin{Proofc}{Proof of Proposition \ref{prop:LlogL}}
  The proof can be directly adapted from the proof of Theorem
  \ref{thm:convexite}. We observe that due to CFL condition (\ref{eqn:CFLexistence}),
  $1-\mu_{i+\frac{1}{2}}^{\alpha,n+1}\in(0,2)$. Using Lemma
  \ref{lem:ftilde}, we obtain the analogue of (\ref{eq:ineglog}): 
  \begin{equation*}
    \sum_{\alpha=1}^{d}{\sum_{i\in\Z}{\tilde{f}(\theta_{i+\frac{1}{2}}^{\alpha,n+1})}} \leq
    \sum_{\alpha=1}^{d}{\sum_{i\in\Z}{\tilde{f}(\theta_{i+\frac{1}{2}}^{\alpha,n})}}
    +
    \sum_{\alpha=1}^{d}{\sum_{i\in\Z}{g(\theta_{i+\frac{1}{2}}^{\alpha,n+1},1-\mu_{i+\frac{1}{2}}^{\alpha,n+1})\ln(1-\mu_{i+\frac{1}{2}}^{\alpha,n+1})}} 
  \end{equation*}
  As $g$ is non-negative, and $\ln(1-\mu)\leq -\mu$ for all $\mu<1$,
  \begin{equation*}
    \sum_{\alpha=1}^{d}{\sum_{i\in\Z}{\tilde{f}(\theta_{i+\frac{1}{2}}^{\alpha,n+1})}} \leq
    \sum_{\alpha=1}^{d}{\sum_{i\in\Z}{\tilde{f}(\theta_{i+\frac{1}{2}}^{\alpha,n})}}
    -
    \sum_{\alpha=1}^{d}{\sum_{i\in\Z}{g(\theta_{i+\frac{1}{2}}^{\alpha,n+1},1-\mu_{i+\frac{1}{2}}^{\alpha,n+1})\mu_{i+\frac{1}{2}}^{\alpha,n+1}}} 
  \end{equation*}
  Using Lemma \ref{lem:ftilde}, we have
  $|\theta-g(\theta,\gamma)|\leq C_2$, for all $\theta\geq 0$ and
  $\gamma\in(0,2)$ with $\gamma_m = 2$. Using equation (\ref{eq:muthetapositive}), 
  \begin{equation*}
    \sum_{\alpha=1}^{d}{\sum_{i\in\Z}{\tilde{f}(\theta_{i+\frac{1}{2}}^{\alpha,n+1})}} \leq
    \sum_{\alpha=1}^{d}{\sum_{i\in\Z}{\tilde{f}(\theta_{i+\frac{1}{2}}^{\alpha,n})}}
    +
    C_2\sum_{\alpha=1}^{d}{\sum_{i\in\Z}{|\mu_{i+\frac{1}{2}}^{\alpha,n+1}|}} 
  \end{equation*}
  As $\mu_{i+\frac{1}{2}}^{\alpha,n+1}=\frac{\Delta t}{\Delta
    x}(\lambda_{i+1}^{\alpha,n+1}-\lambda_i^{\alpha,n+1})$, we get
  \begin{equation*}
    \sum_{\alpha=1}^{d}{\sum_{i\in\Z}{|\mu_{i+\frac{1}{2}}^{\alpha,n+1}|}}\leq
    d\frac{\Delta t}{\Delta x}\Lip(\vect{\lambda})TV(\vect{u}^{n+1})
  \end{equation*}
  Using Lemma \ref{lem:TVD}, we deduce
  \begin{equation*}
    \sum_{\alpha=1}^{d}{\sum_{i\in\Z}{\tilde{f}(\theta_{i+\frac{1}{2}}^{\alpha,n+1})}} \leq
    \sum_{\alpha=1}^{d}{\sum_{i\in\Z}{\tilde{f}(\theta_{i+\frac{1}{2}}^{\alpha,n})}}
    +
    C_2d\frac{\Delta t}{\Delta x}\Lip(\vect{\lambda})TV(\vect{u}^{0}).
  \end{equation*}
\end{Proofc}

\section{Convergence}
\label{sec:convergence}

%%%%%%%%%%%%%%%%%%%%%%%%%%%%%%%%%%%%%%%%%%%%%%%%%%%%%%%%%%%%%%%%%%%%%%%%%%
\subsection{Preliminaries}

We recall the following result (see Lemma 3.2 in \cite{EM1})
\begin{lemma}{\bf ($L\log L$ estimate)}\label{lem::13}\\
If $w \in L^1(\R)$ is a non negative function, then $\int_{\R} \tilde{f}(w) < +\infty$ if and only if 
$w\in L\log L(\R)$. Moreover we have the following estimates
\begin{equation}\label{eq::14}
\int_{\R} \tilde{f}(w) \le 1 + |w|_{L\log L(\R)} + |w|_{L^1(\R)} \ln \left(1+ |w|_{L\log L(\R)}\right)
\end{equation}
\begin{equation}\label{eq::15}
|w|_{L\log L(\R)} \le 1 + |w|_{L^1(\R)}\ln (1+e^2)  + \int_{\R} \tilde{f}(w)
\end{equation}
\end{lemma}

\begin{remark}\label{rem::40}{\bf (Idea of the proof of Lemma \ref{lem::13})}\\
Recall that the proof follows from the following inequalities (for $\mu\in (0,1]$ and $w\ge 0$)
$$\tilde{f}(w)\le w\ln (e + \mu w) + w |\ln \mu| \quad \mbox{with}\quad 1/\mu = 1+ ||w||_{L\log L}$$
and
$$w\ln (e + w) \le 1 + w \ln (1+e^2)+ \tilde{f}(w)$$
We can easily check the following result:
\end{remark}

\begin{lemma}{\bf (Trivial estimate)}\label{lem::15}\\
For any $\gamma\ge 1$ and $\theta \ge 0$, we have
\begin{equation}\label{eq::16}
\tilde{f}(\gamma \theta) \le \gamma \tilde{f}(\theta) + \theta \tilde{f}(\gamma)
\end{equation}
\end{lemma}

We recall the following result (see Lemma 4.3 in \cite{EM1})
\begin{lemma}{\bf (Modulus of continuity)}\label{lem::19}\\
Let $T>0$. Assume that $v\in L^\infty((0,+\infty)\times \R)$ such that
$$|v_x|_{L^\infty((0,T);L\log L (\R))} + |v_t|_{L^\infty((0,T);L\log L (\R))} \le C_1$$
Then for all $\delta,h \ge 0$, and all $(t,x)\in (0,T-h)\times \R$, we have
$$|v(t+h,x+\delta)-v(t,x)|\le 6C_1 \left(\frac{1}{\ln (1+ \frac1{h})}+\frac{1}{\ln (1+ \frac1{\delta})}\right)$$
\end{lemma}

We will state a convergence result in the linear case which will be used later to establish the tame estimate (\ref{eq:tameest}) (see Step 4 of the proof of Theorem \ref{th::2}).
We consider a scalar function $v$ solution of a linear transport equation
\begin{equation}\label{eq::21}
v_t + \lambda^0 v_x =0 \quad \mbox{on}\quad (0,+\infty)\times \R
\end{equation}
where $\lambda^0$ is a real constant and with initial data
\begin{equation}\label{eq::22}
v(0,\cdot)=v_0
\end{equation}
We then consider a solution $v^n$ of an upwind scheme
\begin{equation}\label{eq::23}
\frac{v^{n+1}_{i}-v^{n}_{i}}{\Delta t} 
- \left(\lambda^\varepsilon \right)_- \left(\frac{v^{n}_{i+1}-v^{n}_{i}}{\Delta x}\right)
+ \left(\lambda^\varepsilon\right)_+ \left(\frac{v^{n}_{i}-v^{n}_{i-1}}{\Delta x}\right)=0,\quad \quad \mbox{for}\quad i\in\Z,\quad  n\ge 0
\end{equation}
where $\lambda^\varepsilon$ is a real constant with initial data
\begin{equation}\label{eq::24}
v^{0}_i=v^\varepsilon_i
\end{equation}

\begin{proposition}{\bf (Convergence for the linear scheme)}\label{pro::20}\\
We consider a solution $v$ of (\ref{eq::21})-(\ref{eq::22}) with $v_0\in BUC(\R)$ (the space of bounded and uniformly continuous functions).
We set $\varepsilon=(\Delta t, \Delta x)$ and 
consider the solution $v^n$ to the scheme (\ref{eq::23})-(\ref{eq::24})
for the CFL condition 
$$\frac{\Delta x}{\Delta t}\ge |\lambda^\varepsilon|$$
We set $t_n = n\Delta t$, $x_i=i\Delta x$ and assume that
$$|\lambda^\varepsilon-\lambda^0|\to 0 \quad \mbox{and}\quad \sup_{i\in \Z}|v^\varepsilon_i - v_0(x_i)| \to 0 
\quad \mbox{as}\quad \varepsilon \to 0$$
Then for any compact set $K\subset [0,+\infty)\times \R$, we have with $v^\varepsilon(t_n,x_i)=v^n_i$
$$|v^\varepsilon-v|_{L^\infty(K\cap ( (\Delta t \N) \times (\Delta x \Z))} \to 0 \quad \mbox{as}\quad \varepsilon \to 0$$
\end{proposition}

\begin{Proofc}{Proof of Proposition \ref{pro::20}}
For a viscosity solution $v$, this is an easy adaptation of the general convergence result of Barles, Souganidis \cite{BS}.
It is also easy to check that the limit of the scheme (or directly that $v$) is also a solution in the sense of distributions.
\end{Proofc}

\begin{proposition}{\bf (Weak-$*$ compactness)}\label{pro::30}\\
We consider a sequence of functions $\theta^\varepsilon$ satisfying for some $T>0$
$$|\theta^\varepsilon(t,\cdot)|_{L^1(\R)}+\int_{\R} \tilde{f}(\theta^\varepsilon(t,\cdot)) \le M_T \quad \mbox{for a.e.}\quad t\in (0,T)$$
with $M_T$ a constant independent of $\varepsilon$. Then there exists a function $\theta$ and a constant $C_T = C(M_T)$
such that
\begin{equation}\label{eq::32}
|\theta(t,\cdot)|_{L^1(\R)}+\int_{\R} \tilde{f}(\theta(t,\cdot)) \le C_T \quad \mbox{for a.e.}\quad t\in (0,T)
\end{equation}
such that for any function $\varphi \in C_c((0,+\infty)\times \R)$ (space of continuous functions with compact support), we have
\begin{equation}\label{eq::31}
\int_{(0,+\infty)\times \R} \theta^\varepsilon \ \varphi \to \int_{(0,+\infty)\times \R} \theta \ \varphi
\quad \mbox{as}\quad \varepsilon\to 0.
\end{equation}
\end{proposition}

\begin{Proofc}{First proof of Proposition \ref{pro::30}}
We follow here the lines of the proofs given in \cite{E}.
We consider $\varphi$ with support in $(0,T)\times I$ with $I$ a bounded interval.
We recall that $L\log L(I)$ is defined as $L\log L(\R)$ with $\R$ replaced by the interval $I$.
It is known that $L\log L(I)$ is the dual of $E_{exp}(I) \subset L^\infty(I)$
(see Thm 8.16, 8.18, 8.20 in Adams \cite{A}).
Therefore 
$L^\infty((0,T);L\log L(I))$ is the dual of $L^1((0,T);E_{exp}(I))$ (see Thm 1.4.19 page 17 in Cazenave, Haraux \cite{CH}).
Moreover $L^1((0,T);E_{exp}(I)) \subset L^1((0,T);L^\infty(I))$.
From (\ref{eq::15}), we deduce that
$$|\theta^\varepsilon|_{L^\infty((0,T);L\log L(I))}\le C_{T,I}$$
By general weak-$*$ compactness (see Brezis \cite{Br}), we deduce that for a subsequence, 
there exists a limit $\theta$ (which a priori depends on the compact $[0,T]\times I$, 
but can be chosen independent a posteriori by a classical diagonal extraction argument) such that
(\ref{eq::31}) holds.
Finally, (\ref{eq::32}) follows from Lemma \ref{lem::13}.
\end{Proofc}

\begin{Proofc}{Second proof of Proposition \ref{pro::30}}
We follow the lines of the proofs given in \cite{EM1}.
We recall that from (\ref{eq::15}), we have
$$|\theta^\varepsilon|_{L^\infty((0,T);L\log L(I))}\le C_{T,I}$$
and then using the analogue of Lemma \ref{lem::13} on 
$$A:=(0,T)\times I$$ 
(see Remark \ref{rem::40} for its justification), we have
$$|\theta^\varepsilon|_{L\log L(A)} \le C_{T,I}'$$
for some new constant $C_{T,I}'>0$.
It is known (see page 234 in Adams \cite{A}), that there is a H\"older inequality for the Orlicz space $L\log L (A)$ (with a constant $C$ independent on $A$):
$$||uv||_{L^1(A)} \le C ||u||_{L\log L(A)}  ||v||_{EXP(A)}$$
with
$$||v||_{EXP(A)} = \inf\left\{\lambda >0, \quad \int_A (e^{\frac{|v|}{\lambda}}-1) \le 1\right\}$$
Applying this to $u=\theta^\varepsilon$ and $v=1$, we get that for any measurable set $B\subset A$
$$||\theta^\varepsilon||_{L^1(B)} \le  \frac{C''}{\ln(1+1/|B|)} \quad \mbox{with}\quad C''=C C_{T,I}'$$
This shows that the sequence $\theta^\varepsilon$ is uniformly integrable on $A$, and we can then apply the Dunford-Pettis theorem (see Brezis \cite{Br}), which shows that $(\theta^\varepsilon)_\varepsilon$ is weakly compact in $L^1(A)$, i.e. for any $\varphi\in L^\infty(A)$, we have
$$\int_A \theta^\varepsilon \ \varphi \to \int_A \theta \ \varphi$$
for some function $\theta\in L^1(A)$. In particular this proves
Proposition \ref{pro::30}.
\end{Proofc}

\subsection{Proof of Theorem \ref{th::2}}

\noindent\begin{Proof}
We define $\tilde{S}^n$ the discrete entropy estimate:
\begin{equation*}
  \tilde{S}^n =
  \sum_{\alpha=1}^d{\sum_{i\in\Z}{\tilde{f}(\theta_{i+\frac{1}{2}}^{n,\alpha})\Delta
    x}}
\end{equation*}
\noindent {\bf Step 1: estimate on $\tilde{S}^0$}\\
Using the convexity of $\tilde{f}$, we have with $x_i=i\Delta x$
$$\tilde{f}(\theta^{0,\alpha}_{i+\frac12}) 
= \tilde{f}\left(\frac{1}{\Delta x}\int_{x_i}^{x_{i+1}} (u^\alpha_0)_x(y)dy\right) 
\le \frac{1}{\Delta x}\int_{x_i}^{x_{i+1}} \tilde{f}((u^\alpha_0)_x(y))dy$$
This implies that
$$\tilde{S}^0 \le \sum_{\alpha=1,...,d} \int_{\R} \tilde{f}((u^\alpha_0)_x(y))dy\le C_0$$
where we have used (\ref{eq::14}) to estimate
$$C_0:=\sum_{\alpha=1,...,d} \left\{ 1 + |(u^\alpha_0)_x|_{L\log L(\R)} + |(u^\alpha_0)_x|_{L^1(\R)} \ln \left(1+ |(u^\alpha_0)_x|_{L\log L(\R)}\right)\right\}$$

\noindent {\bf Step 2: Estimates on the $Q^1$ extension $\vect{u}^\varepsilon$}\\
We set $x_i=i\Delta x$ and $t_n = n\Delta t$. Now for $\varepsilon=(\Delta t,\Delta x)$,
we define the $Q^1$ extension of the function defined on the grid, for any $(t,x)\in [t_n,t_{n+1}]\times [x_i,x_{i+1}]$, by 
\begin{multline}\label{eq:Q1ext}
\vect{u}^\varepsilon(t,x) = \left(\frac{t-t_n}{\Delta t}\right)\left\{ \left(\frac{x-x_i}{\Delta x}\right)\vect{u}_{i+1}^{n+1}+ \left(1-\frac{x-x_i}{\Delta x}\right)\vect{u}_{i}^{n+1}\right\}\\
+ \left(1-\frac{t-t_n}{\Delta t}\right)\left\{ \left(\frac{x-x_i}{\Delta x}\right)\vect{u}_{i+1}^{n}+ \left(1-\frac{x-x_i}{\Delta x}\right)\vect{u}_{i}^{n}\right\}
\end{multline}
\noindent {\bf Step 2.1: Estimate on $\vect{u}^\varepsilon_x$}\\
We have for $(t,x)\in [t_n,t_{n+1}]\times(x_i,x_{i+1})$
\begin{equation}\label{eq::34}
\vect{u}^\varepsilon_x(t,x) = \displaystyle \left(\frac{t-t_n}{\Delta t}\right) \vect{\theta}^{n+1}_{i+\frac12} + 
\left(1-\frac{t-t_n}{\Delta t}\right) \vect{\theta}^{n}_{i+\frac12}
\end{equation}
and then using the convexity of $\tilde{f}$
$$\tilde{f}(u^{\varepsilon,\alpha}_x)\le \left(\frac{t-t_n}{\Delta t}\right) \tilde{f}(\theta^{n+1,\alpha}_{i+\frac12}) + 
\left(1-\frac{t-t_n}{\Delta t}\right) \tilde{f}(\theta^{n,\alpha}_{i+\frac12})$$
and then for $t\in [t_n,t_{n+1}]$, we get
\begin{equation}\label{eq::18}
\sum_{\alpha=1,...,d}\int_{\R} \tilde{f}(u^{\varepsilon,\alpha}_x) 
\le \left(\frac{t-t_n}{\Delta t}\right)\tilde{S}^{n+1}+ \left(1-\frac{t-t_n}{\Delta t}\right)\tilde{S}^n   
\le \tilde{S}^0 + C t
\end{equation}
where we have used Proposition \ref{prop:LlogL} for the last inequality.\\
\noindent {\bf Step 2.2: Estimate on $\vect{u}^\varepsilon_t$}\\
Let us define
$$\vect{\tau}^{n+\frac12}_i=\frac{\vect{u}^{n+1}_i-\vect{u}^{n}_i}{\Delta t}$$
We have for $(t,x)\in(t_n,t_{n+1})\times[x_i,x_{i+1}]$
\begin{equation}\label{eq::27}
\left|\vect{u}^\varepsilon_t\right|=\left|\left(\frac{x-x_i}{\Delta x}\right)\vect{\tau}_{i+1}^{n+\frac12}
+ \left(1-\frac{x-x_i}{\Delta x}\right)\vect{\tau}_{i}^{n+\frac12}\right|\leq\left(\frac{x-x_i}{\Delta x}\right)\left|\vect{\tau}_{i+1}^{n+\frac12}\right|
+ \left(1-\frac{x-x_i}{\Delta x}\right)\left|\vect{\tau}_{i}^{n+\frac12}\right|
\end{equation}
Then using the monotonicity and convexity of $\tilde{f}$, we get
$$\tilde{f}(\left|u^{\varepsilon,\alpha}_t\right|)\le 
\left(\frac{x-x_i}{\Delta x}\right)\tilde{f}\left(\left|\tau_{i+1}^{n+\frac12,\alpha}\right|\right)
+ \left(1-\frac{x-x_i}{\Delta x}\right)\tilde{f}\left(\left|\tau_{i}^{n+\frac12,\alpha}\right|\right)\le \tilde{f}\left(\left|\tau_{i+1}^{n+\frac12,\alpha}\right|\right)+\tilde{f}\left(\left|\tau_{i}^{n+\frac12,\alpha}\right|\right) $$
We recall that from the scheme we have with $\lambda^{n+1}_i=\lambda(u^{n+1}_i)$
\begin{equation}\label{eq::28}
\tau^{n+\frac12,\alpha}_i = -(\lambda^{n+1,\alpha}_i)_+\theta^{n,\alpha}_{i-\frac12} 
+ (\lambda^{n+1,\alpha}_i)_-\theta^{n,\alpha}_{i+\frac12}
\end{equation}
and also recall the bound (Theorem \ref{thrm::1} shows that $\vect{u}_i^{n+1}\in\mathcal{U}$)
$$|\lambda^{n+1,\alpha}_i|\le M \quad \mbox{with}\quad M=\max\left(\Lambda^{\alpha},1\right)$$
Therefore applying (\ref{eq::16}), and using the monotonicity of $\tilde{f}$, we get
$$\tilde{f}\left(\left|\tau^{n+\frac12,\alpha}_i\right|\right)\le \sum_{\pm} \left\{M\tilde{f}(\theta^{n,\alpha}_{i\pm \frac12}) + \theta^{n,\alpha}_{i\pm \frac12} \tilde{f}(M)\right\}$$
This implies for $t\in[t_n,t_{n+1}]$
$$\int_{\R}\sum_{\alpha=1,...,d} \tilde{f}(\left|u^{\varepsilon,\alpha}_t\right|)\le
4M\tilde{S}^n + 4 \tilde{f}(M) \sum_{\alpha=1,...,d}\sum_{i\in\Z} \theta^{n,\alpha}_{i+ \frac12} \Delta x$$
We also recall the bound ($\vect{u}^n_i\in\mathcal{U}$)
\begin{equation}\label{eq::29}
|\vect{u}^n_i|\le M_0/(2d)\quad \mbox{with}\quad M_0 = 2d\sum_{\alpha=1}^{d}{\max\left(|m^{\alpha}|,|M^{\alpha}|\right)}
\end{equation}
which implies (using the monotonicity in $x$ of $u^\varepsilon$)
\begin{equation}\label{eq::19}
\sum_{\alpha=1,...,d}
|u^{\varepsilon,\alpha}_x|_{L^\infty((0,T);L^1(\R))}\le M_0
\end{equation}
and
\begin{equation}\label{eq::17}
\int_{\R}\sum_{\alpha=1,...,d} \tilde{f}(u^{\varepsilon,\alpha}_t)\le
4M(\tilde{S}^0 +Ct)+ 4 \tilde{f}(M) M_0
\end{equation}
Moreover, using (\ref{eq::27}), (\ref{eq::28}) and (\ref{eq::19}) we deduce that
\begin{equation}\label{eq:estimateut}
  \sum_{\alpha=1,\dots,d}{|u_t^{\varepsilon,\alpha}|_{L^\infty((0,T),L^1(\R))}}\leq MM_0.
\end{equation}
\noindent {\bf Step 3: Extraction of a convergent subsequence of $\vect{u}^\varepsilon$}\\
From (\ref{eq::19}), (\ref{eq:estimateut}), (\ref{eq::18}), (\ref{eq::17}) and the bound on $\tilde{S}^0$ given in step 1, we see that for any $T>0$, we get the existence of a constant $C_T$ such that
$$\sum_{\alpha=1,...,d}\left\{|u^{\varepsilon,\alpha}_x|_{L^\infty((0,T);L\log L (\R))} + |u^{\varepsilon,\alpha}_t|_{L^\infty((0,T);L\log L (\R))}\right\} \le C_T$$
where we have used (\ref{eq::15}) to estimate the $L\log L$ norm with moreover (\ref{eq::19}).
We also notice that
$$\sum_{\alpha=1,...,d}|u^{\varepsilon,\alpha}| \le M_0/2$$
We can then apply Lemma \ref{lem::19} to get that for any $(t,x)\in (0,T-h)\times \R$, we have
$$\sum_{\alpha=1,...,d}|u^{\varepsilon,\alpha}(t+h,x+\delta)-u^{\varepsilon,\alpha}(t,x)|\le 6C_T \left(\frac{1}{\ln (1+ \frac1{h})}+\frac{1}{\ln (1+ \frac1{\delta})}\right)$$
Therefore by Ascoli-Arzela theorem, we can extract a subsequence (still denoted by $\vect{u}^\varepsilon$) which converges to a limit function $u$ on every compact set $K$ of $[0,+\infty)\times \R$. In particular, we see that the limit function $\vect{u}$ satisfies the initial condition:
$$\vect{u}(0,\cdot)=\vect{u}_0$$ 
Moreover the limit $\vect{u}$ still satisfies
\begin{equation*}\label{eq::26}
\sum_{\alpha=1,...,d}|u^{\alpha}(t+h,x+\delta)-u^{\alpha}(t,x)|\le 6C_T \left(\frac{1}{\ln (1+ \frac1{h})}+\frac{1}{\ln (1+ \frac1{\delta})}\right)
\end{equation*}
\noindent {\bf Step 4: Tame estimate for $\vect{u}$}\\
We want to prove (\ref{eq:tameest}). To this end, we consider a big compact $K$ such that the set
$${\mathcal T}:=\left\{x\ge a + \gamma (t-\tau)\right\}\cap \left\{x\le b - \gamma (t-\tau)\right\}\cap\left\{t\geq\tau\right\}$$
is in the interior of $K$. For any $\varepsilon=(\Delta t, \Delta x)$, we consider 
$i_0\in\Z$ and $N\in\N$ such that
$$[x_{i_0-(N-2)},x_{i_0+(N-2)}]\subset [a,b] \subset [x_{i_0-N},x_{i_0+N}]$$
We consider $n_0\in\N$ and $k\in\N\setminus\{0\}$ such that for $\tau_h=\tau + h$ we have
$$\tau\in [t_{n_0},t_{n_0+1})\quad \mbox{and}\quad  
\tau_h\in [t_{n_0+k}, t_{n_0+k+1}]$$
We recall from (\ref{eq::9}) that we have
\begin{equation*}\label{eq::25}
\frac{1}{k\Delta t} \sum_{\alpha=1,...,d}\sum_{i\in I_{N-k}(i_0)}|u^{n_0+k,\alpha}_i-v^{n_0+k,\alpha}_i| 
\Delta x \le 2\Lip(\vect{\lambda}) \left(TV[\vect{u}^{n_0};I_{N}(i_0)]\right)^2
\end{equation*}
We recall that from Step 3, we have for $(t,x)\in [t_{n_0+k},t_{n_0+k+1}]\times [x_i,x_{i+1}]$
$$\displaystyle \sum_{i=1,...,d}|u^{n_0+k,\alpha}_i - u^\alpha(t,x)|\le 6C_T \left(\frac{1}{\ln (1+ \frac1{\Delta t})}+\frac{1}{\ln (1+ \frac1{\Delta x})}\right) + \sum_{i=1,...,d}|u^{\varepsilon,\alpha} - u^\alpha|_{L^\infty(K\cap ((\Delta t\N)\times (\Delta x \Z))}$$
Using Proposition \ref{pro::20}, this implies in particular that as $\varepsilon\to (0,0)$
$$\frac{1}{k\Delta t} \sum_{\alpha=1,...,d}\sum_{i\in I_{N-k}(i_0)}|u^{n_0+k,\alpha}_i-v^{n_0+k,\alpha}_i| 
\Delta x  \to \frac{1}{h}\sum_{\alpha=1,...,d} \int_{a+\gamma h}^{b-\gamma h} |u^\alpha(\tau +h,x)-v^\alpha(\tau+h,x)|dx$$
with (at least for a subsequence)
$$k\Delta t \to h,\quad k\Delta x \to \gamma h$$
where $\gamma$ can be choosen  bounded in order to satisfy the CFL conditions
(notice that $\gamma$ is also bounded from below, also because of the CFL conditions).
On the other hand, we have
$$TV[\vect{u}^{n_0};I_{N}(i_0)]= \sum_{\alpha=1,...,d} |u^{n_0,\alpha}_{i_0+N+1}-u^{n_0,\alpha}_{i_0-N}|$$
and for the same reasons as previously, we get in particular that
$$TV[\vect{u}^{n_0};I_{N}(i_0)] \to \sum_{\alpha=1,...,d} |u^{\alpha}(\tau,b)-u^{\alpha}(\tau,a)| = TV[\vect{u}(\tau,\cdot);(a,b)]$$
Finally we get
$$\frac{1}{h}\sum_{\alpha=1,...,d} \int_{a+\gamma h}^{b-\gamma h} |u^\alpha(\tau +h,x)-v^\alpha(\tau+h,x)|dx \le 2\Lip(\vect{\lambda}) (TV[\vect{u}(\tau,\cdot);(a,b)])^2$$
which implies (\ref{eq:tameest}). \\
\noindent {\bf Step 5: Passing to the limit in the PDE}\\
\noindent {\bf Step 5.1: Preliminaries}\\
From (\ref{eq:Q1ext}) and (\ref{eq::28}), we have for $(t,x)\in (t_n,t_{n+1})\times (x_i,x_{i+1})$
with $\vect{\lambda}=\vect{\lambda}(\vect{u}(t,x))$ and $\displaystyle a_x = \left(\frac{x-x_i}{\Delta x}\right)$, $\displaystyle b_x=\left(1-\frac{x-x_i}{\Delta x}\right)$
\begin{align}
u^{\varepsilon,\alpha}_t =& \displaystyle  a_x \left\{
-(\lambda^{n+1,\alpha}_{i+1})_+\theta^{n,\alpha}_{i+1-\frac12} 
+ (\lambda^{n+1,\alpha}_{i+1})_-\theta^{n,\alpha}_{i+1+\frac12}\right\} \notag\\
&+ \displaystyle b_x \left\{-(\lambda^{n+1,\alpha}_i)_+\theta^{n,\alpha}_{i-\frac12} 
+ (\lambda^{n+1,\alpha}_i)_-\theta^{n,\alpha}_{i+\frac12}\right\}\notag\\
=&- \lambda_+^\alpha \left\{a_x \theta^{n,\alpha}_{i+\frac12}+ b_x\theta^{n,\alpha}_{i-\frac12}\right\}\notag\\
&+ \lambda_-^\alpha \left\{a_x \theta^{n,\alpha}_{i+\frac32}+ b_x\theta^{n,\alpha}_{i+\frac12}\right\} + e^{\varepsilon,\alpha}(t,x) \label{eq::33}
\end{align}
with
\begin{multline*}
e^{\varepsilon,\alpha}(t,x) = \displaystyle  a_x \left\{
-\left[(\lambda^{n+1,\alpha}_{i+1})_+-\lambda_+^\alpha\right]\theta^{n,\alpha}_{i+1-\frac12} 
+ \left[(\lambda^{n+1,\alpha}_{i+1})_--\lambda_-^\alpha\right]\theta^{n,\alpha}_{i+1+\frac12}\right\}\\
+ \displaystyle b_x \left\{-\left[(\lambda^{n+1,\alpha}_i)_+-\lambda_+^\alpha\right]\theta^{n,\alpha}_{i-\frac12} 
+ \left[(\lambda^{n+1,\alpha}_i)_--\lambda_-^\alpha\right]\theta^{n,\alpha}_{i+\frac12}\right\}
\end{multline*}
In particular, for any test function $\varphi$ with compact support in $K:=[0,T]\times \overline{B_R(0)}$, we have
$$\sum_{i=1,...,\alpha}\left|\int_{[0,+\infty)\times \R} \varphi \ e^{\varepsilon,\alpha}\right|\le 4|\varphi|_\infty T M_0 
\sup_{(\tau,y)\in K}\left(\sup_{|t_{n+1}-\tau|\le \Delta t,\quad |x_i-y|\le \Delta x}
|\vect{\lambda}(\vect{u}^{n+1}_i) - \vect{\lambda}(\vect{u}(\tau,y))|\right)$$
where we have used (\ref{eq::29}). From the uniform convergence of $u^\varepsilon$ on compact sets, we deduce 
in particular that
$$e^{\varepsilon,\alpha} \to 0 \quad \mbox{in}\quad {\mathcal D}'((0,+\infty)\times \R)$$
\noindent {\bf Step 5.2: Introduction of $\vect{\theta}^\varepsilon$}\\
We define the function $\vect{\theta}^\varepsilon$ as
$$\vect{\theta}^\varepsilon(t,x)= \vect{\theta}^n_{i+\frac{1}{2}} \quad \mbox{for}\quad (t,x)\in [t_n,t_{n+1})\times [x_i,x_{i+1})$$
From Proposition \ref{pro::30}, we know that there exists a limit $\vect{\theta}$ such that for any test function $\vect{\varphi}$ (smooth with compact support in $(0,T)\times I$), we have
$$\int_{(0,+\infty)\times \R} \vect{\theta}^\varepsilon \cdot \vect{\varphi} \to \int_{(0,+\infty)\times \R} \vect{\theta} \cdot \vect{\varphi}$$
From (\ref{eq::33}), we also have with $a_x= \displaystyle \frac{x}{\Delta x}- \left\lfloor \frac{x}{\Delta x}\right\rfloor$, $b_x = 1-a_x$
$$u_t^{\varepsilon,\alpha}-e^{\varepsilon,\alpha} 
= -\lambda_+^\alpha \left\{a_x \theta^{\varepsilon,\alpha}+ b_x\theta^{\varepsilon,\alpha}(\cdot,\cdot-\Delta x)\right\}
+ \lambda_-^\alpha \left\{a_x \theta^{\varepsilon,\alpha}(\cdot,\cdot+\Delta x)+ b_x\theta^{\varepsilon,\alpha}\right\}$$
Then
$$A^{\varepsilon,\alpha} = \int_{(0,+\infty)\times \R} (u_t^{\varepsilon,\alpha}-e^{\varepsilon,\alpha})\ \varphi$$
can be computed as follows
$$A^{\varepsilon,\alpha} = \int_{(0,+\infty)\times \R} \theta^{\varepsilon,\alpha} 
\left\{-a_x(\lambda_+^\alpha\varphi) -b_x(\lambda_+^\alpha\varphi)(\cdot,\cdot+\Delta x)  
+a_x(\lambda_-^\alpha\varphi)(\cdot,\cdot-\Delta x)  + b_x(\lambda_-^\alpha\varphi)   \right\}$$
Let us define
$$B^{\varepsilon,\alpha} = \int_{(0,+\infty)\times \R} \theta^{\varepsilon,\alpha} 
\left\{-a_x(\lambda_+^\alpha\varphi) -b_x(\lambda_+^\alpha\varphi)  
+a_x(\lambda_-^\alpha\varphi)  + b_x(\lambda_-^\alpha\varphi)   \right\} = \int_{(0,+\infty)\times \R} -\lambda^\alpha \varphi\ \theta^{\varepsilon,\alpha}$$
Then we have
$$|A^{\varepsilon,\alpha}-B^{\varepsilon,\alpha}|\le \left(\sup_{\pm} 
|(\lambda_\pm^\alpha \varphi)(\cdot,\cdot +\Delta x)- (\lambda_\pm^\alpha \varphi)|_{L^\infty((0,T)\times \R)}\right)\ \int_{(0,T)\times \R} |\theta^{\varepsilon,\alpha}| \to 0$$
On the other hand, we have 
$$B^{\varepsilon,\alpha} \to \int_{(0,+\infty)\times \R} -\lambda^\alpha \varphi\ \theta^{\alpha}$$
This finally shows that
\begin{equation}\label{eq::35}
\forall \alpha \in \{1,\dots,d\}, \, u^\alpha_t + \lambda^\alpha \theta^\alpha=0 \quad \mbox{in}\quad {\mathcal D}'((0,+\infty)\times \R)
\end{equation}
\noindent {\bf Step 5.3: Consequence}\\
Starting from (\ref{eq::34}), we deduce similarly (as in Step 5.2) that
$$u^{\varepsilon,\alpha}_x \to \theta^\alpha \quad \mbox{in}\quad {\mathcal D}'((0,+\infty)\times \R)$$
Therefore
$$\vect{\theta}=\vect{u}_x$$
and from (\ref{eq::35}), we deduce that
$$\forall \alpha \in \{1,\dots,d\}, \, u^\alpha_t + \lambda^\alpha u^\alpha_x = 0 \quad \mbox{in}\quad {\mathcal D}'((0,+\infty)\times \R)$$
with $\vect{u}_x \in L^\infty_{loc}([0,+\infty);L\log L(\R))^d$.\\
\noindent {\bf Step 6: Convergence of the whole sequence when the limit is unique}\\
When we have moreover condition (\ref{eq::8}) for strictly hyperbolic systems,
we know that the solution $\vect{u}$ is unique (among continuous vanishing viscosity solutions).
Therefore, the whole sequence $\vect{u}^\varepsilon$ converges locally uniformly to its unique limit $\vect{u}$.
This ends the proof of the Theorem.
\end{Proof}

%%%%%%%%%%%%%%%%%%%%%%%%%%%%%%%%%%%%%%%%%%%%%%%%%%%%%%%%%%%%%%%%%%%%%%%%
%%%%%%%%%%%%%%%%%%%%%%%%%%%%%%%%%%%%%%%%%%%%%%%%%%%%%%%%%%%%%%%%%%%%%%%%
\noindent {\bf Aknowledgements}\\
RM would like to thank A. El Hajj for useful discussions.

\bibliographystyle{plain}
\bibliography{Article}

\end{document}